\documentclass[9pt,twocolumn,twoside]{pnas-new}

\usepackage{bm}
\usepackage{graphicx}
\usepackage{amsmath, amssymb, amsfonts}
\usepackage{cleveref}
\usepackage{enumitem}
\usepackage{appendix}
\usepackage{xcolor}
\usepackage[shortcuts]{extdash}

\templatetype{pnasresearcharticle} 

\title{Questionnaires to PDEs: From Disorganized Data  \\
to Emergent Generative Dynamic Models}

\author[a]{David W. Sroczynski}
\author[b]{Felix P. Kemeth} 
\author[c]{Ronald R. Coifman}
\author[b,d,1]{Ioannis G. Kevrekidis}

\affil[a]{Department of Chemical and Biological Engineering, Princeton University, Princeton, NJ 08544, USA}
\affil[b]{Department of Chemical and Biomolecular Engineering, Johns Hopkins University, Baltimore, Maryland 21218, USA}
\affil[c]{Department of Mathematics, Yale University, New Haven, CT 06520, USA}
\affil[d]{Department of Applied Mathematics and Statistics, Johns Hopkins University, Baltimore, Maryland 21218, USA}

\leadauthor{Sroczynski} 

\significancestatement{
The extraction of evolutionary partial differential equation (PDE) models from data is predicated on knowing the right {\em independent} variables: space and time coordinates. 
The framework of this paper extends recent developments in machine-learning assisted extraction of PDEs from data to cases where the true time, space, and parameter value instances at which data were obtained are disorganized, and thus unknown. We only know labels for the spatial, temporal and parameter measurement channels. 
Using a Questionnaire metric, the scrambled tensor measurements are organized and smoothly embedded in an intrinsic, {\em emergent} physical space/ time / parameter space, so that generative dynamic models for the data (here, in the form of parameter dependent PDEs) can be learned.
}

\authorcontributions{IGK designed research. DWS, FPK, and IGK performed research, analzed data, and wrote the paper. RRC contributed analytic tools and analyzed data.}
\authordeclaration{The authors declare no conflict of interest.}
\correspondingauthor{\textsuperscript{1}To whom correspondence should be addressed. E-mail: yannisk@jhu.edu}


\begin{abstract}
    Starting with sets of disorganized observations of spatially varying and temporally evolving systems, obtained at different (also disorganized) sets of parameters, 
    we demonstrate the data-driven derivation of parameter dependent, evolutionary partial differential equation (PDE) models capable of generating the data. 
    This tensor type of data is reminiscent of shuffled (multi-dimensional) puzzle tiles.
    The {\em independent variables} for the evolution equations (their ``space'' and ``time'') as well as their effective parameters are all ``emergent'', i.e., determined in a data-driven way from our disorganized observations of behavior in them.
    We use a diffusion map based ``questionnaire'' approach to build a parametrization of our emergent space/time/parameter space for the data.
    This approach iteratively processes the data by successively observing them on the ``space'', the ``time'' and the ``parameter'' axes of a tensor.
    Once the data are organized, we use machine learning (here, neural networks) to approximate the operators governing the evolution equations in this emergent space.
    Our illustrative example is based on a previously developed vertex-plus-signaling model of \textit{Drosophila} embryonic development. 
    This allows us to discuss features of the process like symmetry breaking, translational invariance, and autonomousness of the emergent PDE model, as well as its interpretability.  
\end{abstract}

\dates{This manuscript was compiled on \today}
\doi{\url{www.pnas.org/cgi/doi/10.1073/pnas.XXXXXXXXXX}}

\begin{document}

\maketitle
\thispagestyle{firststyle}
\ifthenelse{\boolean{shortarticle}}{\ifthenelse{\boolean{singlecolumn}}{\abscontentformatted}{\abscontent}}{}

\dropcap{D}ata science and machine learning daily expand the set of data-driven tools in the mathematical modeler's toolkit.
This toolkit enables, among other tasks, the extraction of data-driven dynamic models 
capable of predicting the evolution of a system's response as a function of initial conditions and (possibly) external parameters.
The input to this process is a (rich enough) set of time series (or image series, movies)
of experimental observations of the system we wish to model.

\begin{figure*}
    \includegraphics[width=17.8cm]{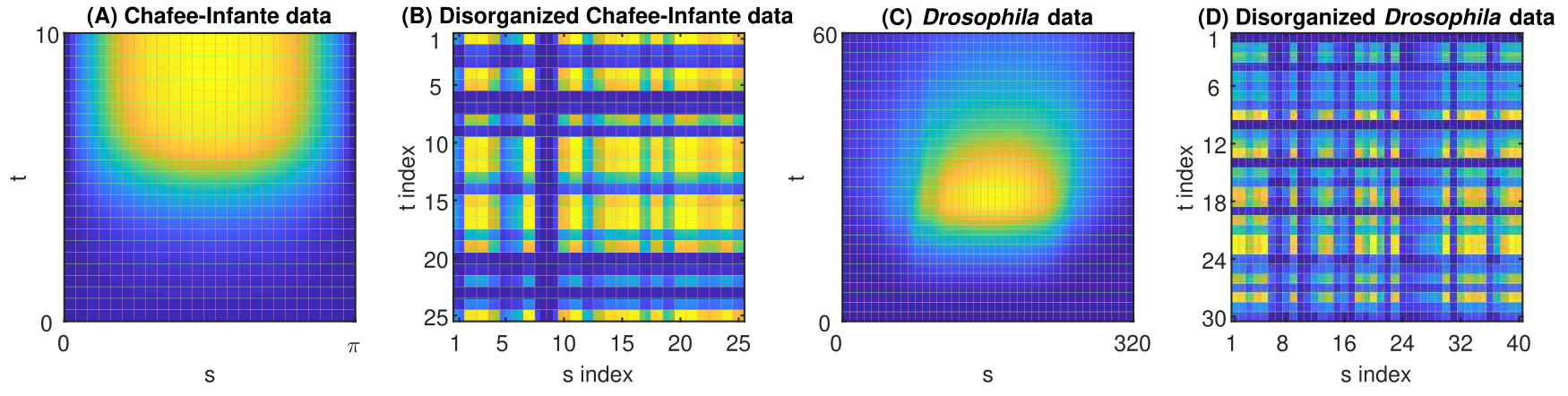}
    \caption{\label{fig:CI_scramble_plus_embryo} Example space-time plots. (A) A space-time plot for the evolution of the \mbox{Chafee-Infante} partial differential equation ($u_t=u-u^3+\nu u_{xx}$) for $\nu=0.16$, with Dirichlet boundary conditions $u=0$ at $x={0,1}$ and initial condition $u=0.1$ for all x (excepting the boundary). (B) Chafee-Infante space-time data disorganized in space and time. (C) A space-time plot for chemical signal intensity data in a \textit{Drosophila} embryonic development model. (D) The \textit{Drosophila} data disorganized in space and time.}
\end{figure*}

\begin{SCfigure*}
    \includegraphics[width=11.4cm]{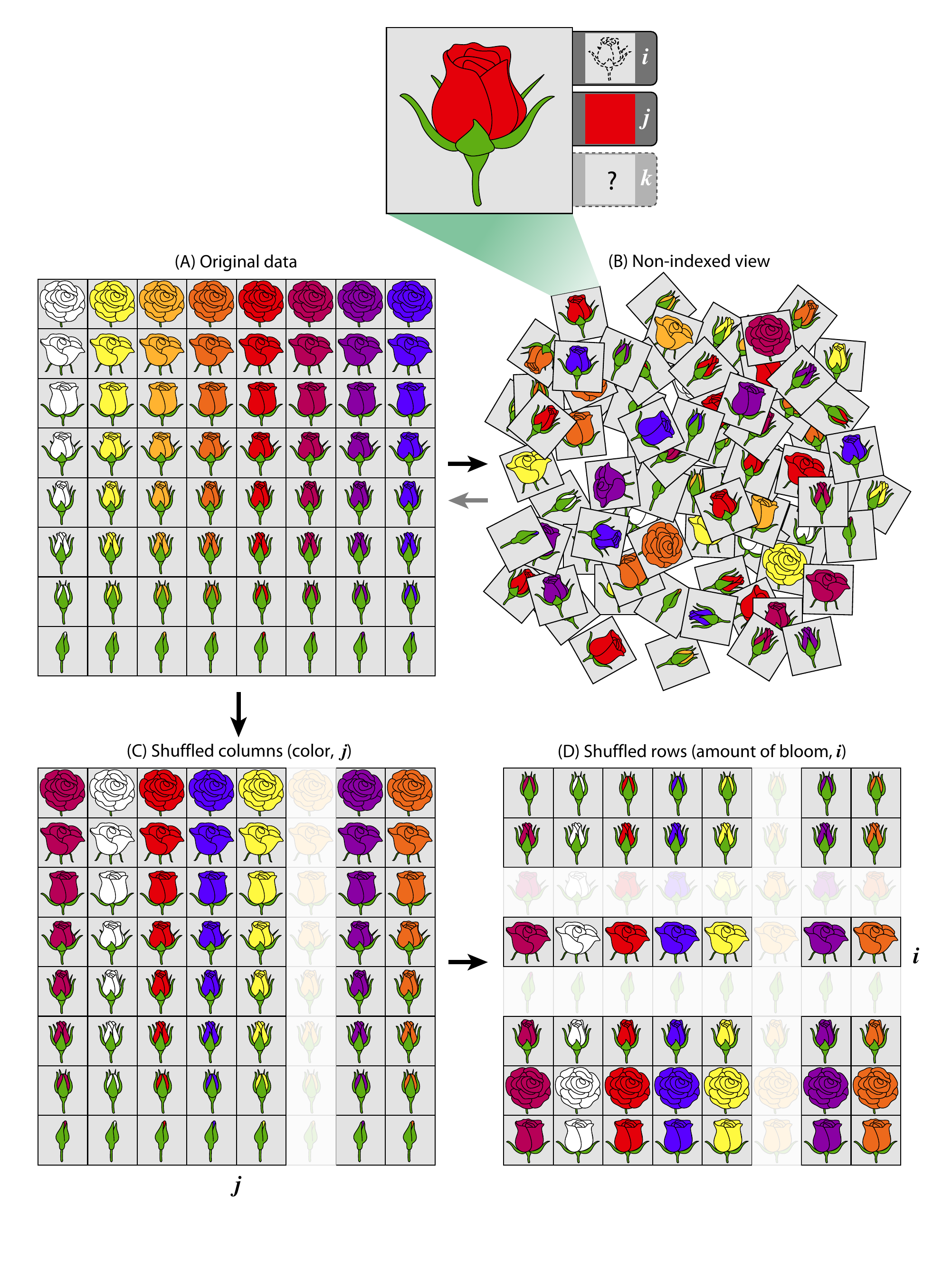}
    \caption{\label{fig:rose_scramble} A caricature of data shredding and scrambling represented by roses in different colors and stages of bloom. \textit{(A)} The full data field. Each grid box represents a potentially multi-dimensional measurement. The rows represent the sequential age of the rose (i.e., time), while the columns represent the color of the rose (i.e., a graded parameter). \textit{(C)} The column indices $j$ are shuffled, so the colors are out of order. We may also have some data missing (represented by the faded column) which will have to be interpolated. \textit{(D)} Now the row indices $i$ are shuffled, resulting in the times being non-sequential. We may again have some data missing (represented by the faded rows) which will have to be interpolated. However, note that roses of the same age are still in the same row, and roses of the same color are still in the same column. This version of shuffling (which is relevant for this paper) is more ordered than the fully jumbled tiles in \textit{(B)}. This is because each rose image is ``tagged'' \textit{(Top)} with a row index $i$ and a column index $j$ (and possibly as we show later in the paper, with a third index $k$); roses of the same color index are still together, and roses of the same age index are still together, even if we don't know the color that corresponds to a specific $j$ index, nor the time that corresponds to a specific $i$ index.} 
\end{SCfigure*}

A simple illustration is seen in Fig. \ref{fig:CI_scramble_plus_embryo}(A): a space-time plot of the evolution of the \mbox{Chafee-Infante} equation, a scalar, one-dimensional partial differential equation (PDE), with fixed boundary conditions $u=0$ at $x=0,1$ and initial condition $u=0.1$ for all x (excepting the boundary), from time $t=0$ to time $t=10$.

Given this ``movie'' in the form of the space-time field $u(x,t)$ we can obtain (at every x and every t) a set of measurements:  $u, u_t, u_{tt}, u_x, u_{xx}, u_{xxx}$, etc. 
If we have some reason to believe that the dynamics can be modeled in the form of a partial differential equation of the form $u_t = \mathcal{L}\left(u, u_{xx}\right)$, then each point of the movie gives us a point in the $u_t, u, u_{xx}$ space. It is clear that, with these data, the operator $\mathcal{L}$ can be approximated (fitted) as a function of ``just'' $u$ and $u_{xx}$, and any off-the-shelf neural network or Gaussian Process Regression software  can be used to ``learn the right-hand-side of the PDE''. 
Cautionary notes abound: notice, for example,  how much we have already assumed, even without a formula: that it is a first order PDE in time; that it is translationally invariant (its law does not explicitly depend on $x$) and autonomous (its law does not explicitly depend on $t$); that the right-hand-side operator only depends on $u$ and its second order spatial derivatives; that no other variable is necessary to predict the evolution of $u$; that noise can be ignored...
Yet the fact remains: in this data-driven sense, operators can be approximated, and ODEs and PDEs can be ``learned'' from data through, say, neural networks; this has been known and practiced for decades\cite{Krischer1993,rico-martinez92_discr_vs,Gonzalez-Garcia1998} (and is experiencing an explosive rebirth in the current literature\cite{Brunton2020,Lu2021,Karniadakis2021,Rudy2017,zhang_ma_2020,Raissi2018}).
These learned models do not need to be completely ``black box'' agnostic: physical knowledge can be included in ``hardwiring'' parts of the operator that are assumed accurately known, or learning the ``calibration'' of parts of the operator that are assumed only partially/qualitatively known. The term ``gray box identification'' is used for such algorithms\cite{Rico-Martinez1994,Lovelett2020}.

Figure \ref{fig:CI_scramble_plus_embryo}(C) is a qualitatively similar computational space-time movie: it arises in modeling the evolution in time of a chemical signal from a vertex-plus-chemistry model of \textit{Drosophila} egg evolution (proposed in Ref. \citenum{HocevarBrezavscek2012} and summarized in Appendix \ref{app:drosophila}). It only records a particular observation of the evolution (along a portion of the one-dimensional ``backbone'' of the equatorial slice of the egg). Without any of the myriad details, the point is that one could try and ``learn'' an evolution PDE for this spatiotemporal signal from the data.

Such a data-driven model can only be guaranteed, upon successful training, to be a compact summary of the data it was trained on: it can reproduce them (it can regenerate the data in the same way that a PDE solver can produce a solution for a well-posed problem). How well it can generalize (extrapolate at other initial conditions, other boundary conditions, other parameters) or whether it can assist physical understanding is, of course, another story that only starts after the small initial success of creating a compact ``generator'' of the training data. 
We will take all this ``compact data generator'' technology for granted, and use it in our work here.

A first sketch of the problem we want to solve is outlined in Figure \ref{fig:CI_scramble_plus_embryo}(B): 
Figure \ref{fig:CI_scramble_plus_embryo}(A) has been turned into a ``shredded and shuffled'' puzzle.
Measurements (pixels, puzzle tiles) are obtained at $N_s$ locations in space ($s_i, i = 1,...,N_s$); yet the instruments are placed at random space locations, so that while the $i^{th}$ instrument always measures at the same spatial location, we do not know where this location
is in physical space.
Similarly, the $N_s$ instruments are triggered to record at the same $N_t$ instances in time ($t_j, j = 1,...,N_t$); yet the labels $j$ of the temporal measurements are also random (not sequential with time); so, all $N_s$ spatial channels report at time instance $j$, but we
do not know at what physical time these simultaneous measurements were taken.
%
We thus have a list of ``spatial channels'' and a a list of ``temporal channels'' that index our observational tiles (and in what follows, we will make the puzzle three-dimensional: we will add ``parameter channels,'' performing
$N_p$ different experiments with disorganized labels $k = 1, N_p$; this will turn our ``tiles'' into ``voxels''). 

So: we know what measurements were obtained simultaneously in time (from their index $j$); which come from the same
experiment (from their index $k$); and which come from the same spatial location (from their index $i$); but 
we do not know what the actual physical time corresponds to the index $j$, which particular space location corresponds to the index $i$, and for what parameter values the measurements at parameter index $k$ were obtained; our tensor data are ``trebly disorganized.'' 

The caricature in Fig. \ref{fig:rose_scramble} reiterates what we mean here by ``shuffled'' or ``disorganized''. The ``solved puzzle'' in Fig. \ref{fig:rose_scramble}(A) represents the full, organized data field we want to be able to reproduce. Each grid box represents a (potentially multi-dimensional) data point. The rows represent sequential time (i.e., the degree at which the rose has bloomed), while the columns represent a gradually changing parameter (i.e., the color of the rose). As we move from Fig. \ref{fig:rose_scramble}(A) to (C), the column indices $j$ are shuffled, and the organization of gradual color change is ``lost''. Moving to (D), the row indices $i$ are shuffled, and the organization in sequential time is ``lost''. Note (full disclosure!) that this scrambling is still more ordered than the random ``pile of tiles'' shown in Fig. \ref{fig:rose_scramble}(B): the set of measurements from a given row remain together (now in a different row), and likewise measurements from a given column remain together (now in a different column). This is possible because each measurement is ``tagged'' with indices $i$ for its time channel and $j$ for its parameter channel (i.e., the label of the experiment in which it was observed). In other examples, either $i$ or $j$ might represent a spatial index rather than time or parameters. As we will show later, one can include a third index $k$, in which case we have a three-dimensional data tensor including space, time, and parameters viewpoints all at once. It is also possible for the data to be incomplete (i.e. we may not have access to observations of data {\em everywhere} in the domain, but rather at a selection of points in space/time/parameter space). Areas with missing observations appear faded in Fig. \ref{fig:rose_scramble}(C/D); we will show how to ``fill them in'' as necessary. Alternative illustrations of the shuffling can be found in Appendix \ref{app:questionnaire}, Fig. \ref{fig:triangle_scramble}.

We want to combine (a) organizing/reconstructing the (multi)puzzles (finding a good way to embed our measurement locations and temporal instances in an ``emergent space-time'' domain, 
or an ``emergent space-time-parameter'' domain) with (b) learning a generative model (an evolution equation) {\em in
this reconstructed, data-driven, ``emergent space-time''} with the machine-learning-assisted techniques mentioned above.

While the literature of reconstructing dynamical systems from data (with known space-times!) is quite rich and fast growing \cite{gonzalez-garcia98_ident_distr_param_system,Thiem2021}, there is also 
extensive literature for the mathematics and algorithms of (instantaneous!) puzzle reconstruction \cite{Huroyan2020,Sholomon2016}.
Our algorithm of choice here for ``solving the puzzle'' (i.e. for the reconstruction of a useful spatial-temporal-parameter embedding from observational ``puzzle voxels'') will be the tensor decomposition ``Questionnaires'' algorithm of Ref. \citenum{Ankenman2014}.
This is illustrated, and mathematically summarized, in Appendix \ref{app:questionnaire}; it was proposed in Ref. \citenum{Ankenman2014}
and we have used it in the past to learn normal forms from dynamical system observations in Ref. \citenum{Yair2017,Sroczynski2018} (see also \cite{Moon2019,Gigante2019}).
Our algorithm for subsequently learning a generative model for the organized data -an evolution equation in the ``emergent space-time-parameter'' domain-  will be deep neural networks (as proposed originally in Refs. \citenum{Hudson1990,Rico-Martinez1992,Rico-Martinez1994,Rico-Martinez1995,Gonzalez-Garcia1998,Anderson1996} and used extensively recently, e.g. in \cite{Arbabi2020,Arbabi2021,Lee2020,Kemeth2020}.

The paper is organized as follows:
We will start with a very brief description of Diffusion Maps, and their use as part of the Questionnaires algorithm (detailed and illustrated in Appendix \ref{app:questionnaire}). We will then briefly introduce the data we will use, which come from a previously proposed/studied vertex-plus-signalling model of \textit{Drosophila} embryonic development, described in detail in Appendix \ref{app:drosophila}. Finally, we will describe and illustrate ``learning the PDE'' in the emergent domain.
When a problem is even mildly nontrivial, interesting twists arise in treating it; here, these twists include (a) slight breaking of a left/right symmetry and (b) the fact that we know that the problem is not spatially translationally invariant: it includes a physically motivated, spatially localized, temporally varying forcing term (the source of the signalling). How these two ``twists'' arise and are dealt with in our data-based scheme is, we believe, of some interest.

We conclude with a few of the (myriad) caveats and shortcomings of the approach.
Even with those, we  argue that the combination of puzzle-solving with nonlinear distributed system identification, and the ability to 
create ``intelligent'' emergent domains in which to learn smooth models, is an important pursuit, extending the tools of modern data-driven modeling.

A final note before starting: why scramble a space-time you already know?. The answer is that we first validate the approach on problems where we know the solution, before it can be applied to data with hidden space-times (think of segments of broken fossils in different earth strata at different locations,
as well as data where no obvious physical space-time exists (e.g. dynamics of networks, e.g power networks or physical neural networks)), but which can be usefully visualized in data-driven space-times (e.g. Ref. \citenum{Mishne2016,Kemeth2018}).

\section{METHODS FOR DATA ANALYSIS}
\label{computational methods}
This section briefly introduces the data organization tools: the {\em Diffusion Maps} manifold learning technique, as well as the {\em informed metric} that iteratively synthesizes information viewed along different axes of the data tensor.

\subsection{Manifold Learning: Diffusion Maps}
The goal of manifold learning is to discover underlying lower-dimensional intrinsic nonlinear structure in high-dimensional data. Diffusion maps\cite{Coifman2006,Lafon2004} accomplishes this by constructing a discrete approximation of the Laplacian operator on the data. When the data are sampled from a low-dimensional manifold, the discrete operator converges to the continuous Laplace-Beltrami operator {\em on the manifold} in the limit of infinite sample points. The discrete operator is constructed by defining a weighted graph on some $N$ sampled data points, where the weight between points $i$ and $j$ is given by
\begin{equation}
\label{dmaps}
w_{i,j}=\exp\left(- \frac{d(\bm{y}_i,\bm{y}_j)^2}{\epsilon^2} \right);
\end{equation}
$d(\bullet,\bullet)$ represents a chosen distance metric (e.g. Euclidean in the ambient space) and $\epsilon$ represents a distance scale below which samples are considered similar. A weight of 1 indicates that two samples are identical, while a weight close to 0 indicates that two samples are very dissimilar. After some normalization, the leading non-harmonic eigenvectors $\{\phi_k\}$ of the kernel matrix, weighted by the corresponding eigenvalues $\{\lambda_k\}$, provide a new coordinate system for embedding/describing the data. Distances in this coordinate system are referred to as {\em diffusion distances}. Eigenvectors which do not contribute to this distance (due to low eigenvalues) can be truncated, and the reduced set of eigenvectors can serve as a proxy for the intrinsic manifold coordinates. More details can be found in the references\cite{Coifman2006,Lafon2004}.

\subsection{Iteratively Informed Geometry and the Questionnaire Metric}
One of the key choices in the implementation of diffusion maps is that of the metric used to compare data points. In many cases, the standard Euclidean norm in ambient space is sufficient, but in certain applications other metrics may warranted. Returning to the caricature in Fig. \ref{fig:rose_scramble}, consider the case where we want to find a jointly smooth embedding for the rose color channels as well as for the blooming stage (age) channels at which we collect data. The questionnaire metric (see Appendix \ref{app:questionnaire} and Ref. \citenum{Yair2017, Sroczynski2018,Mishne2016}) uses (1) the data-driven geometry of the blooming stage indices to inform the distances between the color channel indices; as well as (2) the data-driven geometry of the color channel indices to inform the distances between the blooming stage indices. The procedure iteratively improves the joint metric until convergence. If a third ``viewpoint'' (in addition to color and blooming stage - e.g., possibly type of fertilizer used) is included, we iterate between all three viewpoints.

\section{Illustrative Example}

\subsection{Model Data}

To illustrate our approach, we will use data generated by a vertex-plus-signalling model of a tubular epithelium which approximates ventral furrow formation in early \textit{Drosophila} embryos; the model was first developed in Ref. \citenum{HocevarBrezavscek2012}, see also Ref. \citenum{Polyakov2014,Misra2016,Misra2016a}). The approach combines:
\begin{enumerate}
    \item A 2D mechanical model consisting of a ring of 80 quadrilateral cells. Since each pair of neighboring cells shares two vertices, the state space of the mechanical model is described by the positions of 160 vertices. The vertices are acted on by line tension on the edges, a stiff outer membrane, and energy penalties for deviation in the volume of each individual cell, as well as the central ``yolk'' they collectively envelop.
    \item A chemical signal model for a protein involved in embryonic development. The (temporally varying) chemical intensity of the signal is assumed to be spatially uniform within each cell (the cells are ``well mixed''). There is a source for the signal in certain cells, whose intensity is time-dependent: it grows to a maximum value at the ``stopping time'' $t_s$ before decaying exponentially with first order rate constant $d$. Transport between the cells is proportional to the concentration difference between them, and is characterized by an ``effective diffusivity'' $D_e$.
\end{enumerate}

The overall model overlays the mechanical and chemical components to generate ``videos'' (series of snapshots) in the spirit of experimentally tracking the staining for the relevant protein. For more details, see Appendix \ref{app:drosophila}.
The model contains a number of constitutive parameters; here we will fix several of them, and focus on {\bf three} that we will allow to vary: the effective diffusivity $D_e$, the protein degradation rate constant $d$, and the stopping time $t_s$.

We generated data from this \textit{Drosophila} embryo model for $1000$ distinct parameter settings. For each configuration, $D_e$ and $d$ were generated from independent normal distributions with \{$\mu_{D_e}=0.2,\sigma_{D_e}=0.04$\} and \{$\sigma_d=0.075,\sigma_d=0.005$\}. Settings beyond two standard deviations from the mean were discarded and redrawn. 
The stopping time $t_s$ was then taken to be a prescribed function of $D_e$ and $d$; the point of this is to illustrate that, {\em even though \underline{three parameters are varying}, there is only a \underline{two-parameter family of variations}} - so that our parameter settings lie on a two-dimensional manifold in the three-dimensional parameter space (Fig. \ref{fig:param_embed}(A)). We will thus expect our data-driven approach to recognize that the effective parameter variation is two-dimensional. 

\begin{figure*}
    \includegraphics[width=17.8cm,height=0.6\textheight]{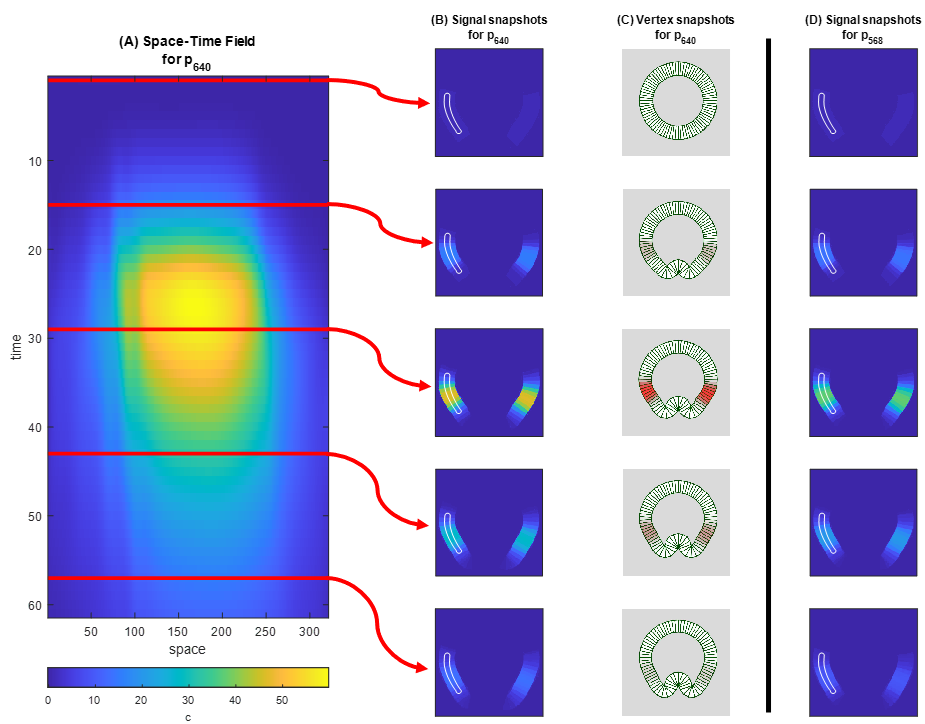}
    \caption{\label{fig:spacetimeVSsnapshot} (A) A space-time {\em Drosophila} embryo chemical signal field for a particular parameter setting ($d=0.1750, D_e=0.0650, t_s=40.17$). Five representative points in time are highlighted in red. (B) Observation snapshots corresponding to those points in time. Data is taken from the area outlined in white. (C) Corresponding snapshots, showing the location of the cell vertices. (D) Signal Snapshots at the same time instances but for a different parameter setting ($d=0.1830, D_e=0.0848, t_s=39.91$).}
\end{figure*}

A representative simulation output is summarized in Fig. \ref{fig:spacetimeVSsnapshot}(C).
From this type of output we generated ``videos'' of just the chemical signal expression (Fig. \ref{fig:spacetimeVSsnapshot}(B)), with no specific morphology information, since the morphology evolves in the same way in all our trials.
%
In order to simplify analysis
, we sample data from a one-dimensional ``backbone curve'', guided by the midpoints of the cell interface edges (see white outlines in Fig. \ref{fig:spacetimeVSsnapshot}(B/D)).

We approximate the data, sampled on $N=360$ points along this curve of cell midpoints,
using bivariate splines. With $T=61$ snapshots for each movie, this results in a $1000 \times 61 \times 360$ data tensor;
each snapshot contains $N=360$ spatial grid points.
Fig. \ref{fig:spacetimeVSsnapshot} shows how a space-time field for a particular parameter setting corresponds to snapshots from the original video;
to illustrate variability, some snapshots from a different parameter setting are also included
in the last column.

\subsection{Embedding Results}

\begin{figure*}
    \includegraphics[width=17.8cm]{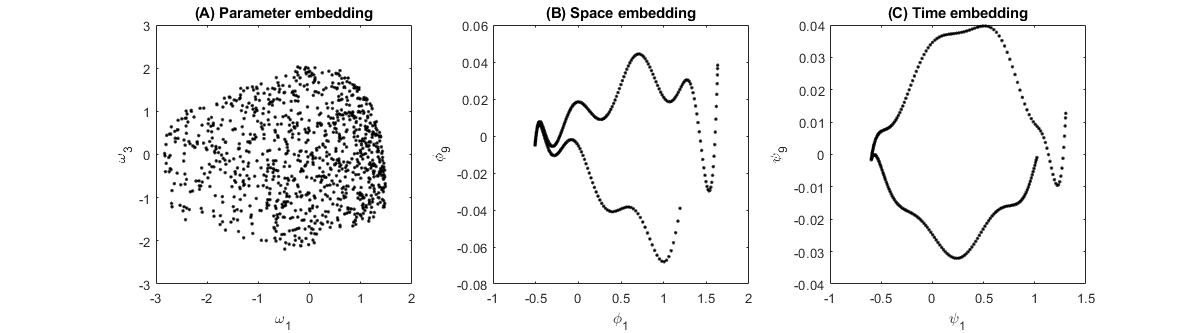}
    \caption{\label{fig:all_embed} The recovered questionnaire embeddings for the (A) parameters, (B) space, and (C) time (see text).}
\end{figure*}

After applying diffusion maps with the questionnaire metric to our \emph{scrambled} \textit{Drosophila} data, a subset of the diffusion map eigenvectors provide an embedding for each axis of our data tensor (parameters, space, and time), which are shown in Fig. \ref{fig:all_embed}. While these embeddings may not be intuitive at first glance, they can (in this example) be rationally interpreted in terms of the underlying system.

\begin{figure*}
    \includegraphics[width=17.8cm]{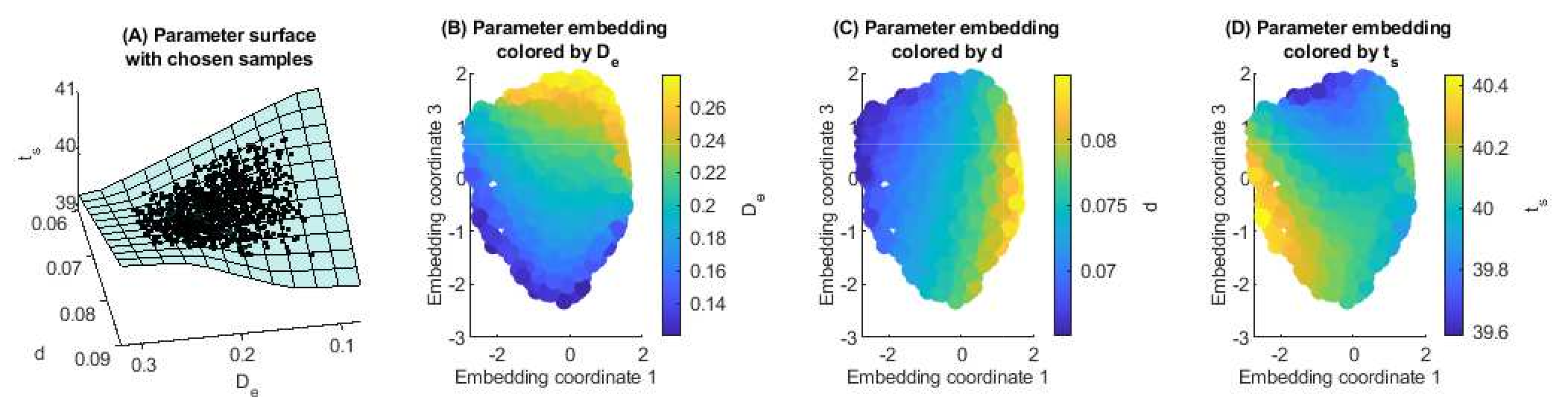}
    \caption{\label{fig:param_embed} (A) Parameter settings of the observations (black dots) in {$D_e$,$d$,$t_s$}-space, with the blue surface showing the 2D manifold those parameters were drawn from.
(B, C, D) Data-driven parameter embedding coordinates colored by $D_e$, $d$, and $t_s$. The embedding is one-to-one with the two-dimensional parameter domain surface.}
\end{figure*}

We begin with the embedding of the parameters, since it is the simplest. Since our parameter settings were sampled from a two-dimensional manifold, we would expect the algorithm to embed the data with only two unique coordinates. In general, taking diffusion map eigenvectors with the highest eigenvalue may not give the most parsimonious embedding, since ``higher harmonics'' of significant coordinates may appear before new unique coordinates. Methods exist, however, to filter out such unnecessary coordinates \cite{Dsilva2018}. In this case, the first and third diffusion map eigenvectors were the only relevant coordinates, with eigenvector 2 being a function (``higher harmonic'') of eigenvector 1. Fig. \ref{fig:param_embed} shows that the two recovered coordinates are bi-Lipschitz with the true parameters, with the first coordinate being mostly a function of $d$ and the other mostly a function of $D_e$. Thus, with no prior assumptions on the nature of the system, we have established the effective two-dimensionality of the parameter sampling, and we have revealed the underlying intrinsic organization geometry of the data in parameter space.

\begin{figure*}
    \includegraphics[width=\textwidth,height=0.31\textheight]{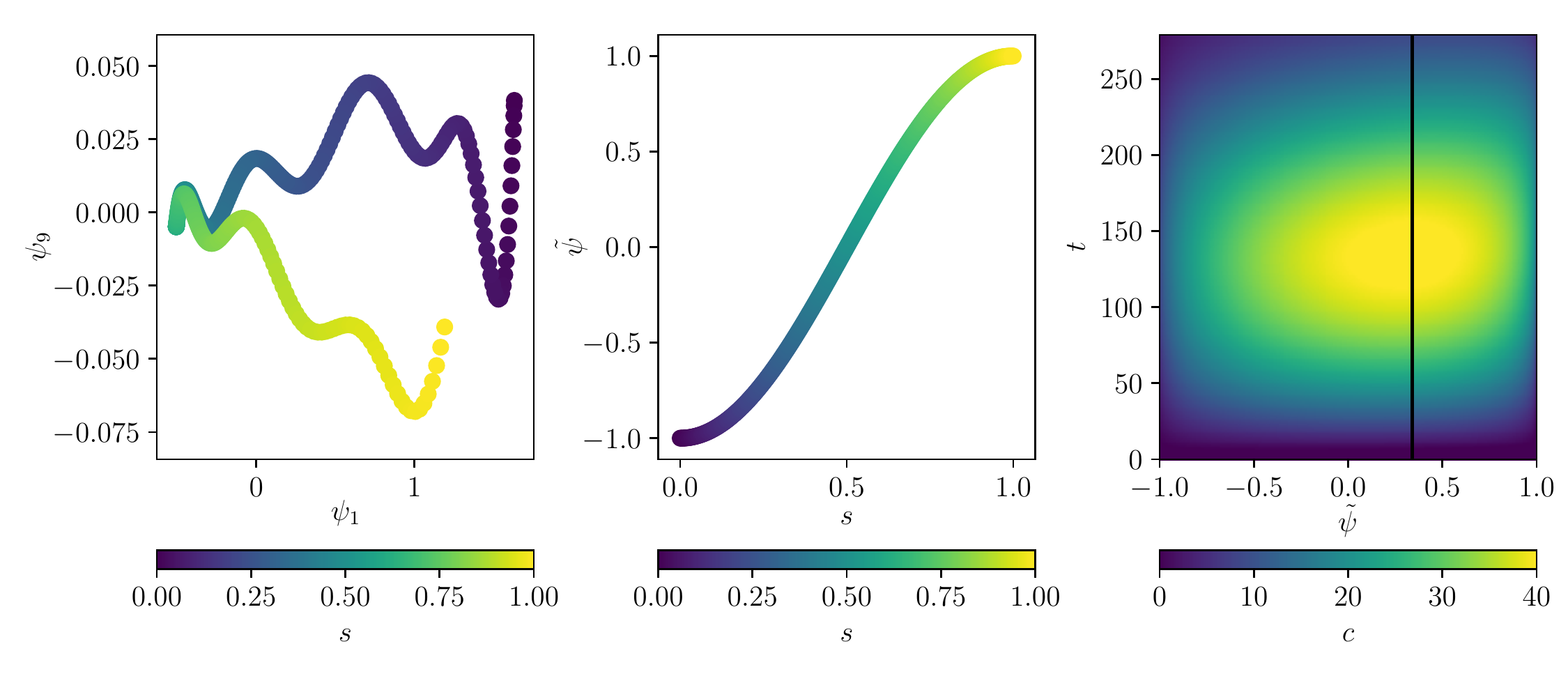}
    \caption{Left: Data-driven space embedding obtained from the questionnaires, colored with the position $s$ along the cell section centerline from which the data was taken.
      Middle: Arc length $\tilde{\psi}$ extracted from the space embedding $\psi$, as a function of the
      position $s$.
      Right: Concentration data parametrized by the extracted embedding arclength $\tilde{\psi}$.
      The black vertical line corresponds to the position of the minimum of $\psi_1$; that is,
      to the left edge of the embedding shown on the left.}
    \label{fig:learn_space}
\end{figure*}

Given that the data is sampled from a 1D curve in space, it is reasonable to expect a 1D embedding for space; yet in this case the algorithm gives a 2D space embedding (Fig.~\ref{fig:all_embed}(B)) which is locally 1D.
Note the ridged ``hairpin''-like shape of this embedding, which can be explained as follows:
For a perfectly circular domain, reaction-diffusion dynamics on both sides of a source cell would give left-right symmetric concentration profiles.
However, the asymmetric (and moving!) source cell locations introduce a symmetry breaking into the system, making
the branches ``above'' and ``below''  the source close {\em but distinguishable}. This becomes clear in the space embedding
(Fig.~\ref{fig:learn_space}(left)),
where $\psi_1$ encodes the distance to the source term and $\psi_9$ the induced symmetry breaking.

We therefore use the arclength along this ``hairpin'' to parametrize the emergent spatial geometry of the data.
An effective coordinate $\tilde{\psi}$ is extracted using diffusion maps on the curve in Fig. \ref{fig:learn_space}(left) with a nearest neighbor similarity measure, 
and shown in Fig.~\ref{fig:learn_space}(middle) as a function of the arclength $s$ along the cell centerline from which data was taken - notice the one-to-one correspondence.
We use this coordinate $\tilde{\psi}$ as the emergent, data-driven space coordinate~\cite{kemeth2020_learning_pdes_in_emergent_coordinates} in which to learn the dynamics of $c(\tilde{\psi}, t)$ (Fig.~\ref{fig:learn_space}(right)).

The embedding of the time samples (Fig. \ref{fig:all_embed}(C)) also has a similar quirk, in that it requires two dimensions. The first embedding coordinate roughly follows the overall chemical concentration, which starts at zero, rises to a maximum near $t_s$, and then fades back asymptotically to zero. Because the numerical experiment is stopped at finite time, the embedding coordinate never reaches its original value, but comes close. There is a need for a second embedding coordinate ($\phi_9$) which captures the fact that the spatial distribution of chemical concentrations is different when the overall level is rising (more tightly focused around the source cells) than when it is falling (more spread out due to diffusion between cells). In this case the diffusion is not that strong, so the difference is slight, which is why this second coordinate ``shows up'' later on in the spectral hierarchy as $\phi_9$. Essentially, the time evolution of this system has been characterized by the algorithm as a skinny loop.

\begin{figure*}
    \includegraphics[width=\textwidth,height=0.31\textheight]{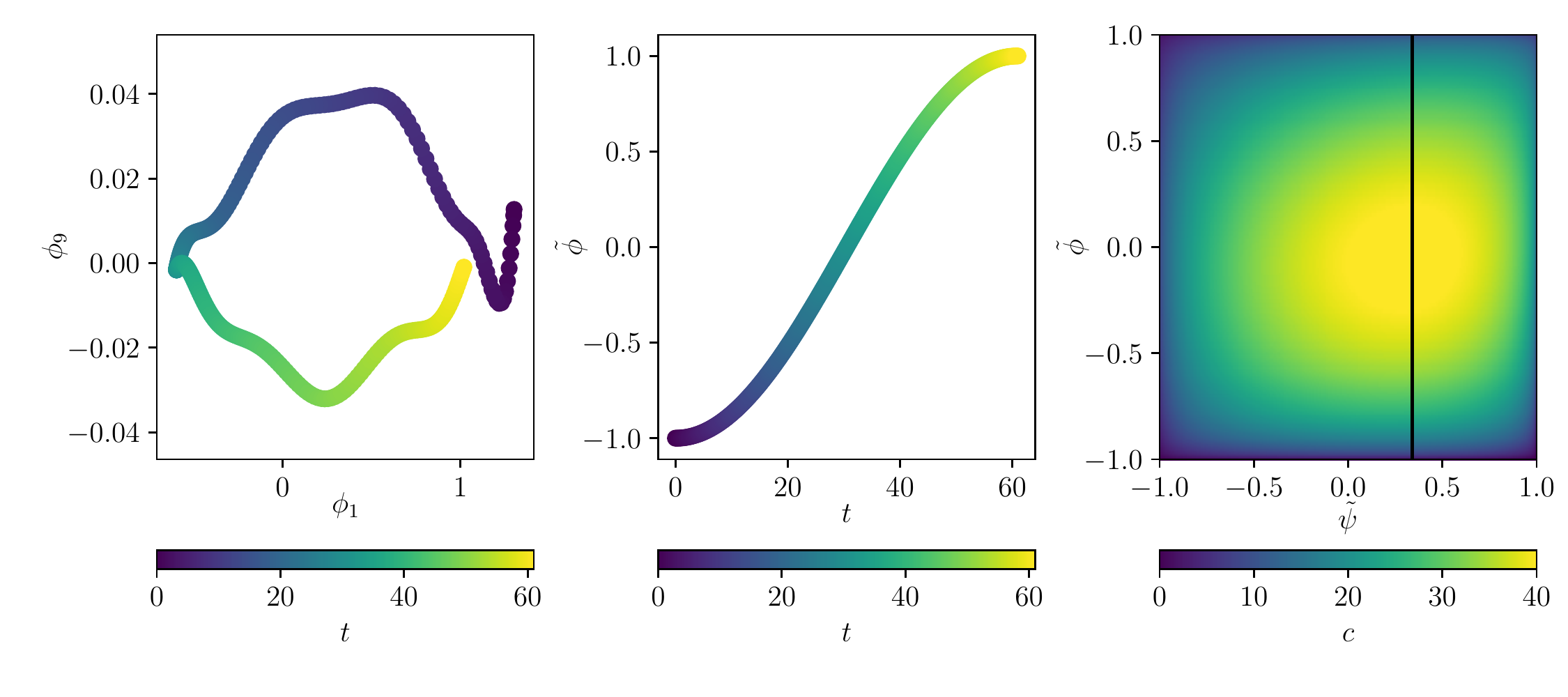}
    \caption{Left: Data-driven time embedding, obtained from the questionnaires, colored with the true time $t$.
      Middle: Embedding arclength $\tilde{\phi}$ extracted from the time embedding $\phi$, as a function of the true time $t$.
      Right: Concentration data parametrized by the extracted arclengths $\tilde{\psi}$ and $\tilde{\phi}$.
      The black vertical line corresponds to the position of the minimum of $\psi_1$ as above.}
    \label{fig:learn_time}
\end{figure*}

Similarly to our approach to the space embedding, we can extract a coordinate $\tilde{\phi}$ from this skinny loop, to obtain a 1D embedding for emergent time (Fig. \ref{fig:learn_time}(middle)). We also show the data in this final ``emergent'' space-time in Fig.\ref{fig:learn_time}(right).


\subsection{Learning the dynamics in the emergent coordinates}

We start the section by approximating the evolution operator in emergent space but still in physical time; we pose a distributed parameter model (a PDE) whose right hand side is approximated 
by a convolutional (in emergent space) neural network~\cite{gonzalez-garcia98_ident_distr_param_system, rico-martinez92_discr_vs}. We approximate the dynamics of the concentration for a single parameter setting through a PDE
\begin{equation}
  \frac{\partial c}{\partial t} = f\left(c, \frac{\partial c}{\partial \tilde{\psi}}, \dots, 
  \frac{\partial^n c}{\partial \tilde{\psi}^n}\right)
\label{eq:model}
\end{equation}
where $f$ is represented by a neural network and $c=c\left(\tilde{\psi}, t\right)$. Here, we use $n=4$ derivatives, estimated using finite differences. We therefore resample the data on an equally spaced grid in $\tilde{\psi}$ using a bivariate spline approximation, and also on $T=1500$ equally-spaced points in actual time.

$f$ is composed of 6 fully connected hidden layers with 126 neurons each, with one output layer containing a single node. Each hidden layer is followed by a Swish activation function~\cite{ramachandran17_swish}.
The model is optimized using the PyTorch framework~\cite{paszke2019pytorch} and a Adam~\cite{kingma2017adam} optimizer with the default hyperparameters, based on the mean squared error of the predicted and true $\frac{\partial c}{\partial t}$.
The true time derivatives are estimated using finite differences in time.
The initial learning rate is set to 0.005, and subsequently halved when the training loss did not decrease for 75 epochs.
The model is trained for a total of 1500 epochs using a batch size of 128. Overfitting was assessed using a held out validation set composed of all spatial points of the respective last ten snapshots.
Note that we do not provide a source term in the model. We therefore learn the PDE in a ``corridor'' around the source location, obtained as described above.
Furthermore, boundary conditions are not in principle available for the learned model. We therefore provide, in lieu of boundary conditions, narrow ``boundary corridors'' informed by the data.
As in Ref.~\cite{kemeth2020_learning_pdes_in_emergent_coordinates}, we regularize the outputs of the learned model using a truncated singular value decomposition. Finally, we integrate an initial $c$ profile using the learned model.

\begin{figure*}
    \includegraphics[width=\textwidth]{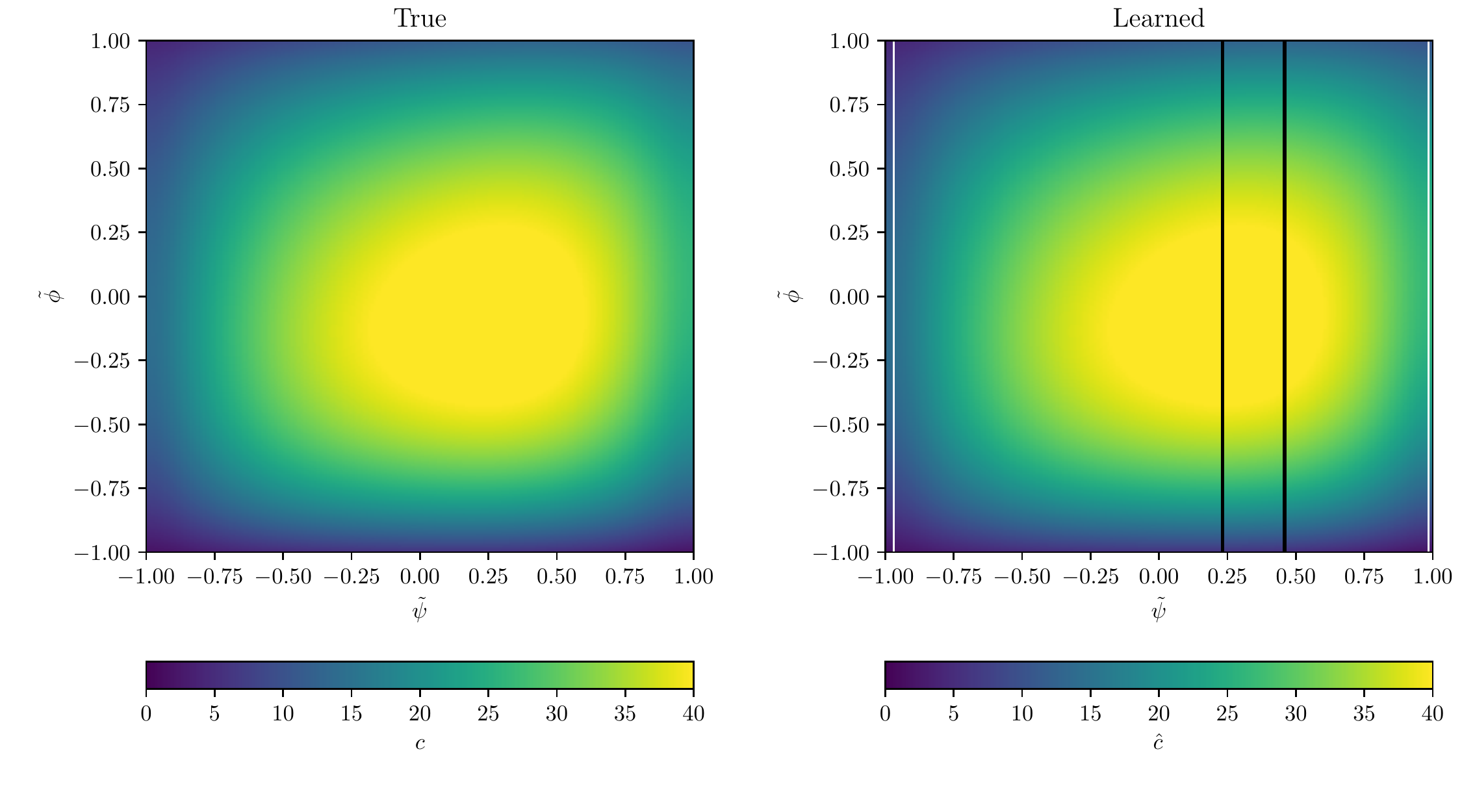}
    \caption{Left: True data as shown in Fig.~\ref{fig:learn_time}(right). 
Right: Integration results using the non-autonomous ML-learned emergent PDE model, ~\eqref{eq:f_w_source}, using the same initial snapshot. Note the visual agreement with the true data,}
    \label{fig:learn_results_emergent_space_and_time}
\end{figure*}

Note that here, we learned the dynamics, cf. Eq.~\eqref{eq:model}, at just a single parameter setting, and took the temporal ordering of the snapshots as known and given. 
However, if the true times are not known, we can instead construct the model to integrate in {\em emergent} time. As discussed above, time gives a hairpin embedding $\phi$ similar to the space embedding, so we use a similar emergent coordinate $\tau=\tilde{\phi}$ (see Fig.~\ref{fig:learn_time}), and learn a model in this emergent time

\begin{equation}
  \frac{\partial c}{\partial \tilde{\phi}} = f\left(c, \frac{\partial c}{\partial \tilde{\psi}}, \dots, \frac{\partial^n c}{\partial \tilde{\psi}^n}\right)
  \label{eq:tau_model}
\end{equation}
with $c=c\left(\tilde{\psi}, \tau\right)$. The results of integrating this model are shown in Fig.~\ref{fig:learn_results_emergent_space_and_time}. 
Data-informed boundary corridors, as well as the source term are provided.

If the source term is not given, the predictions of the learned (translationally invariant) partial differential equation
model are simply wrong. This becomes obvious when using the learned model to predict
the dynamics over the entire domain (only providing boundary conditions at the edges) but not the source information. Alternatively, one can extend this approach {\em to learn the source term in the corridor $I_s$} as well; for further results and discussion of learning a non-autonomous dynamical system, see Appendix \ref{app:non-autonomous}.

Future work will focus on incorporating the emergent parameter coordinates, $\omega$, into the learning process,
\begin{equation}
  \frac{\partial c}{\partial \tau} = f\left(c, \frac{\partial c}{\partial \tilde{\psi}}, \dots, \frac{\partial^n c}{\partial \tilde{\psi}^n}; \boldsymbol{\pi}\right),
\end{equation}
providing a fully data-driven system identification framework.

\section{CONCLUSIONS}

The mathematics underpinning the data-driven solution of puzzles have started, in recent years, to provide increasingly sophisticated puzzle reconstructions, including cases of missing data. Even cases where different parts of the puzzle have been observed through different sensors (so, ``puzzle fusion'') are starting to appear.
Our purpose here was to combine a data organization technology (Questionnaires) with the (machine-learning assisted) construction of generative mathematical models. More specifically, our generative models came in the form of differential equations (here, partial differential equations): models which, given a few initial/boundary conditions, allow us to reconstruct the entire puzzle.
In this sense, what we present can be thought of as a combination of data organization and ``boosted'' data compression: now, with very few data (initial/boundary conditions) and a dynamic
generative law (here, a parabolic PDE) we can reconstruct good approximations of all the missing data, and even sometimes extrapolate successfully.
We should stress again that what we did would be {\em much easier} if we had explicit time/space information; here we had to invent the data-driven, emergent space-time in which the data appear smooth, and where, therefore, a parabolic PDE type model can be postulated. 
This ``boosted'' data compression, in the form of ``very few data plus generative law'' can now be used to interpolate in parameter space or in emergent physical space-time; one may even attempt to extrapolate (up to when singularities will arise).
%

More importantly, the approach naturally allows for ``physics infusion'' -- if one has an informed guess of what the actual independent space variable should be, or of what an approximate closed form of the generative dynamic model could be, this information can be included in the process in a ``gray box'' identification scheme that will {\em calibrate} the partial physical knowledge to the quantitative truth (the data) in the form of a multifidelity calibration problem.

It is important also to note how some crucial assumptions (homogeneity of the emergent space, or autonomousness (homogeneity in emergent time))  shape the entire process; if there is reason to believe they do not hold, then fitting a space-time homogeneous generative model to the data 
will reproduce the (training) data, interpolate between them, {\em and fail miserably upon generalization/extrapolation}. Our approach cannot therefore proceed successfully unless we have sufficient {\em qualitative} information/hypotheses about the problem symmetries, so as to formulate a reasonably well-posed optimization problem.

What about interpretability, and what about understanding? The first is much easier: we can always create data-driven ``mirrors"  (calibration curves, surfaces, hypersurfaces, homeomorphisms, diffeomorphisms, conjugacies, dualities) that map the data-driven predictions to a comparable number of physically interpretable ones. We only need to test that these ``translation engines'' (these diffeomorphisms between data-driven observables and physically interpretable ones) have a well-behaved Jacobian on the data (it is bounded away from zero and from infinity). 
Notice that many equally successful (on the data) physical interpretations may be possible - one can model a lot of gravity data equally well with Newtonian gravity and with general relativity; on such data, we can say that the two theories are ``reconcilable'': they can be calibrated to each other.
Among different successful interpretations (different sets of physically interpretable quantities that are also one-to-one with the data-driven observables), we can even choose the one with the best numerical behavior (best condition number; best Lipschitz constants).

But the crucial question remains:  we can predict, but what did we \underline{understand} about the generating dynamics and about the underlying physical mechanisms? One answer is ``nothing, emphatically nothing!''. Our procedure looks like a crystal ball - we ask questions and obtain answers, but do not have a {\em mechanistic} understanding that can be described {\em parsimoniously} in terms of physically interpretable observables, functions, and operators. Even more - we do not have equations that can be, upon inspection, decomposed into the ``reaction'' part, the ``elastic'' part and the ``diffusion'' part (alternatively, into the 
``gradient'' part, the ``harmonic'' part and the ``rotational'' part, to make an analogy with the Helmholtz-Hodge decomposition). 
Often one resorts here to sparse identification: whatever function fits the data well with ``a few'' common dictionary terms ``ought to be the right physical interpretation''. We understand this, but respectfully disagree. It is (at least a little) presumptuous to expect the truth to be parsimoniously expressible in our everyday favorite dictionary.
That does not mean we do not revere Occam's razor. 
But we believe that, ultimately, understanding might be more in terms of the algorithms that discover the generative relations across the data, and not necessarily, or at least not exclusively, in discovering a few compact algebraic terms that parametrize these relations reasonably accurately. If the latter {\em can} be done, there is plenty of evidence that it sometimes leads to good generalizations, good extrapolations, good ``disentanglements'' of phenomena. 
Yet ultimately, the truth does not {\em have to} appear beautiful in our own favorite current language. 
Maybe it is {\em the language in which truth is beautiful} that we should strive to formulate - the transformation to the space in which the evolution is isospectral, as Peter Lax would say in his Lax Pair formulation, or the space in which ``troubles melt like lemon drops'', as Dorothy would sing in the Wizard of Oz.

\matmethods{Please describe your materials and methods here. This can be more than one paragraph, and may contain subsections and equations as required. Authors should include a statement in the methods section describing how readers will be able to access the data in the paper. 

\subsection*{Subsection for Method}
Example text for subsection.
}


\acknow{The authors are grateful to Prof. S. Shvartsman and Dr. M. Misra for stimulating discussions and for graciously providing their simulation code for {\em Drosophila} epithelial shape transformations. They are also grateful to Dr. O. Yair and Prof. R. Talmon of the Technion for the original implementation of questionnaires used. This work was partially supported by the US AFOSR (FA9550-21-1-0317,FA9550-20-1-0288) and by the DARPA Atlas program (W911NF21C0010).}

\showacknow{}

\bibliography{main}


\pagebreak
\appendix

\begin{figure*}
\title{\textbf{SUPPLEMENTARY MATERIAL}

\textbf{FOR}

\textbf{Questionnaires to PDEs: From Disorganized Data to Emergent Generative Dynamic Models}
}
\\
David W. Sroczynski\textsuperscript{a}, Felix P. Kemeth\textsuperscript{b}, Ronald R. Coifman\textsuperscript{c}, Ioannis G. Kevrekidis\textsuperscript{b,d,1}

\textsuperscript{a}Department of Chemical and Biological Engineering, Princeton University, Princeton, NJ 08544, USA;
\textsuperscript{b}Department of Chemical and Biomolecular Engineering, Johns Hopkins University, Baltimore, Maryland 21218, USA;
\textsuperscript{c}Department of Mathematics, Yale University, New Haven, CT 06520, USA;
\textsuperscript{d}Department of Applied Mathematics and Statistics, Johns Hopkins University, Baltimore, Maryland 21218, USA

\textsuperscript{1}To whom correspondence should be addressed. E-mail: yannisk@jhu.edu
\end{figure*}

\section{Iteratively Informed Geometry and the Questionnaire Metric}
\label{app:questionnaire}

\begin{figure*}
    \includegraphics[width=17.8cm]{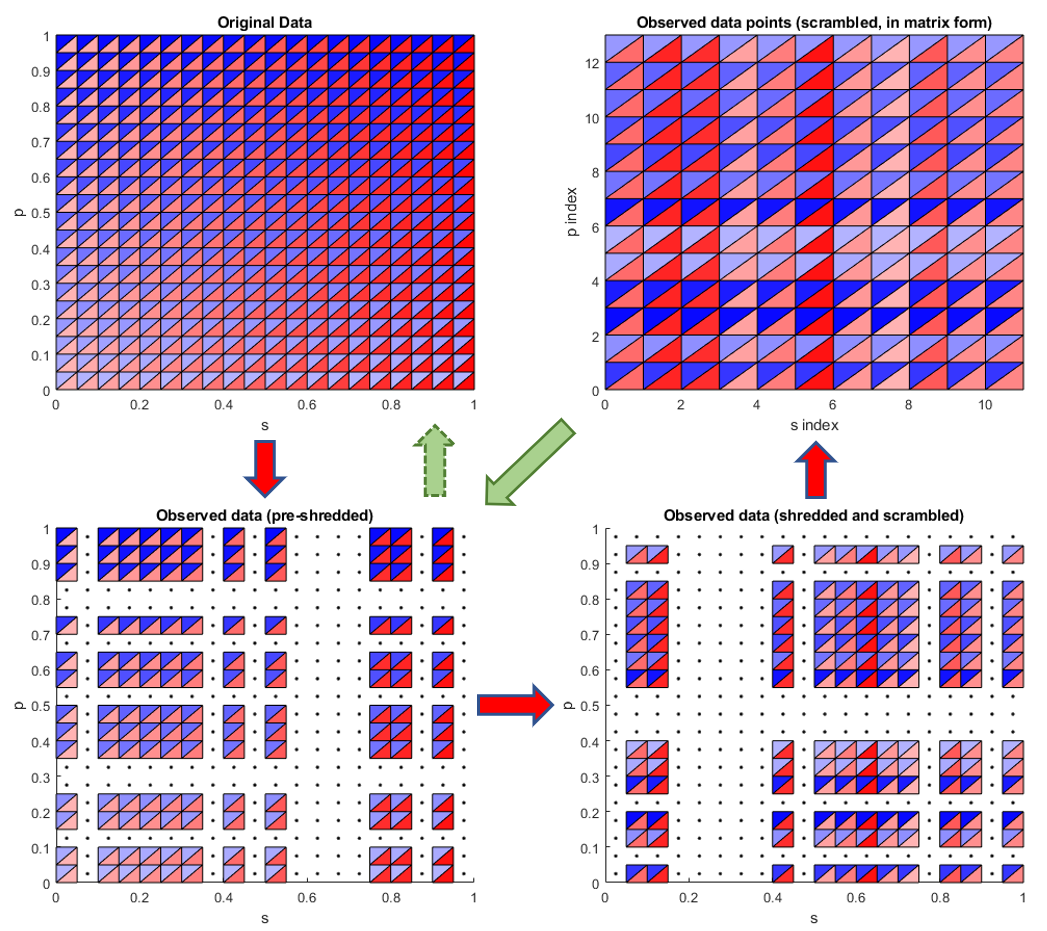}
    \caption{\label{fig:triangle_scramble} Caricature of data shredding and scrambling. \textit{(Top Left)} The full data field in one spatial variable $s$ and one parameter $p$. Note that the red intensity is constant within a column, and the blue intensity is constant within a row. \textit{(Bottom Left)} We only have access to the field at a random sampling of values for $s$ and $p$. \textit{(Bottom Right)} By scrambling the indices for space and parameters, we ``shred'' the data. However, the red intensity is still constant within a column, and the blue intensity is still constant within a row. \textit{(Top right)} After we have removed empty rows and columns, this data is the input to our algorithm. \textit{(Solid green arrow)} We will use diffusion maps with the questionnaire metric to unscramble our data, both reordering them, and placing them correctly in physical space and parameter space (more accurately in new, emergent space and parameter embedding coordinates, that are one to one with physical space and parameter space). \textit{(Dashed green arrow)} Various interpolation methods can be used to fill in the missing areas of the field.}
\end{figure*}


The concept behind an {\em informed metric}, such as the {\em Questionnaire metric} for tensor-type data measured with different {\em types of channels} (the ``axes of the tensor''), is that we can construct better similarity measures between data points by incorporating information about the geometries {\em of the different types of channel measurements} at the points.
Suppose we want to find an embedding space for the parameters that affect the behavior of a dynamical system, based on observed system behavior. 
We do not assume any knowledge of the underlying dimensionality of the right parameter space, nor how the parameters enter into the governing equations. 
Instead, we have data for the system behavior generated from several unknown parameter settings (drawn, say, from some distribution). 
Regardless of how the many parameter settings are generated, for each parameter setting we collect several measurements of the system behavior using several channels. 
%
If the system is described by a spatial field, each channel might measure the field variable at a different location in space. In this scenario, our data comes in the form of an $N_p$ by $N_c$ matrix, where $N_p$ is the number of parameters and $N_c$ is the number of spatial measurement channels. If $\mathbf{y}_{\mathbf{p}_i}$ is the column-vector of $N_c$ measurements for parameter setting $\mathbf{p}_i$, then the $i^{th}$ row of our data matrix is $\mathbf{y}_{\mathbf{p}_i}^T$. However, this is only one “viewpoint” of the matrix; we could similarly denote the $i^{th}$ column of our data matrix by $\mathbf{y}_{\mathbf{c}_i}$, the vector of measurements of channel $\mathbf{c}_i$ at all $N_p$ measured parameter settings. It is important to note that the parameter settings ${p_i}$ are disordered; a given $i$ identifies a particular parameter setting (experiment), but says nothing about how that parameter setting relates to any other.  Figure \ref{fig:triangle_scramble} presents a caricature of the shredding and scrambling of data in ``parameter space'' and ``physical space'' for this case. The description that follows is based on the thesis of Ankenman\cite{Ankenman2014} and previous illustrations of its use can be found in, e.g., Ref. \citenum{Yair2017,Sroczynski2018,Mishne2016}.

When constructing an embedding of the parameter settings, the input to diffusion maps is the set of pairwise distances $d(\mathbf{y}_{\mathbf{p}_i},\mathbf{y}_{\mathbf{p}_j})$ for all $\mathbf{p}_i$,$\mathbf{p}_j$, for some distance metric $d(\bullet,\bullet)$. Common choices for $d(\bullet,\bullet)$ include the Euclidean and the $L_1$ norms. While there are different ways to incorporate the geometry of the {\em spatial} measurement channels into an informed metric, the questionnaire metric does so by appending additional coordinates to the $\mathbf{y}_{\mathbf{p}_i}$ and $\mathbf{y}_{\mathbf{p}_j}$ data vectors of the ``parameter view'' before taking a standard $L_1$ norm; this gives the form
\begin{equation}
    d_{quest}(\mathbf{y}_{\mathbf{p}_i},\mathbf{y}_{\mathbf{p}_j})={||\mathbf{y}_{\mathbf{p}_i}-\mathbf{y}_{\mathbf{p}_j}||}_1+{||\mathbf{F}^p (\mathbf{y}_{\mathbf{p}_i})-\mathbf{F}^p (\mathbf{y}_{\mathbf{p}_j})||}_1,
\end{equation}

\noindent where $\mathbf{F}^p (\mathbf{y}_{\mathbf{p}_i})$ is the vector of additional coordinates appended to $\mathbf{y}_{\mathbf{p}_i}$. The $k^{th}$ element of $\mathbf{F}^p (\mathbf{y}_{\mathbf{p}_i})$ is given by

\begin{equation}
   \mathbf{F}^p (\mathbf{y}_{\mathbf{p}_i})=\langle \mathbf{g}_k,\mathbf{y}_{\mathbf{p}_i} \rangle,    
\end{equation}

\noindent where $\mathbf{g}_k$ is a basis vector {\em defined over the $N_c$ {\bf spatial} measurement channels}. In this form, $\mathbf{F}^p$ is a transform of $\mathbf{y}_{\mathbf{p}_i}$ using a set of basis vectors \{$\mathbf{g}_k$\} based on the spatial measurement channels (in the questionnaire metric, the set of basis vectors will be overcomplete).

To find these basis vectors in a data-driven manner that incorporates the intrinsic {\em spatial channel geometry}, we construct a set of hierarchical clusters of the spatial measurement channels using a bottom-up approach based on some distance metric $d(\mathbf{y}_{\mathbf{c}_i},\mathbf{y}_{\mathbf{c}_j})$ between the spatial channels. At the bottom level, each channel belongs to its own cluster. At each successive level, we form super-clusters by attempting to join each individual cluster to either another individual cluster or an existing super-cluster. We start by joining the two individual clusters with the smallest distance between their respective members (since at the start there are no super-clusters), and we continue choosing the individual cluster with the smallest average distance to another individual cluster or supercluster. We only actually join these clusters if the distance between them is lower than some threshold distance, which grows larger at each successive level. We continue until we have attempted to join each individual cluster once (and only once); after which we move to the next level (at the next level, the superclusters from the lower level are relabeled as individual clusters). At the highest level, there will be one root cluster containing all spatial channels. Note that other clustering methods are also acceptable.
The result of this clustering is a set of clusters \{$\mathbf{I}_k$\}, each containing some subset of the channels. We then define the basis vectors

\begin{equation}
\mathbf{g}_k(\mathbf{c}_i)=
\left\{
    \begin{array}{ll}
        1, & \mathbf{c}_i \in \mathbf{I}_k, \\
        0, & otherwise.
    \end{array}
\right.
\end{equation}

\noindent Since $\mathbf{g}_k$ is an indicator function for whether $\mathbf{c}_i$ is in $\mathbf{I}_k$, $\mathbf{F}_k$ is the sum of values in $\mathbf{y}_{\mathbf{p}_i}$ corresponding to the channels belonging to $\mathbf{I}_k$. {\em It has been shown\cite{Ankenman2014} that defining the questionnaire metric in this fashion is equivalent to an earth-mover's distance between $\mathbf{y}_{\mathbf{p}_i}$ and $\mathbf{y}_{\mathbf{p}_j}$.}
So far, we have not specified the distance metric used in the clustering of the spatial channels; for our first iteration, we simply use an $L_1$ norm. 

We could also define, in a completely analogous way, a {\em questionnaire metric} for the spatial ``view'' of the tensor:
\begin{equation}
    d_{quest}(\mathbf{y}_{\mathbf{c}_i},\mathbf{y}_{\mathbf{c}_j})={||\mathbf{y}_{\mathbf{c}_i}-\mathbf{y}_{\mathbf{c}_j}||}_1+{||\mathbf{F}^c (\mathbf{y}_{\mathbf{c}_i})-\mathbf{F}^c (\mathbf{y}_{\mathbf{c}_j})||}_1.
\end{equation}

\noindent As the parameter view questionnaire metric required clustering of the {\em spatial} measurement channels, the spatial view questionnaire metric requires clustering of the {\em parameter} measurement channels, based on some distance metric. Having, in the first iteration step, constructed our informed metric for the parameter view, {\em we can now use it for augmented parameter vector clustering}. This parameter-view clustering is used to construct additional coordinates for our spatial data vectors $\mathbf{y}_{\mathbf{c}_i}, \mathbf{y}_{\mathbf{c}_j}$; and 
through this augmentation, to construct an {\em informed spatial metric}.
Now the procedure iterates: spatial clustering helps update the informed parameter metric; parameter clustering helps update the informed spatial metric. 
This line of reasoning leads to the following iterative algorithm:

\begin{enumerate}[nosep]
	\item Cluster the spatial channels using an uninformed metric.
	\item Use the uninformed spatial channel clusters to construct an informed metric and to perform clustering for the parameter channels
	\item Repeat until converged:
	\begin{enumerate}[nosep]
	    \item Use the parameter clusters to update the spatial channel metric and spatial channel clustering
	    \item Use the spatial channel clusters to update the parameter channel metric and parameter channel clustering
	\end{enumerate}
\end{enumerate}

This process is typically observed to converge within a few iterations\cite{Yair2017,Sroczynski2018,Mishne2016}.

\subsection{3-D Questionnaire Metric}

Now suppose that we measure each channel not just once, but at a series of $N_t$ time points. Our data matrix is now an $N_p$ by $N_c$ by $N_t$ tensor. Now there are three viewpoints of our data, and we can define a questionnaire metric for each one based on a hierarchical clustering of the other two. 
We define $\mathbf{Y}_{\mathbf{p}_i}$ as the $N_c$ by $N_t$ matrix of all spatial channel measurements at all time channels for parameter $\mathbf{p}_i$. 
In the definition of the questionnaire metric, $\mathbf{y}_{\mathbf{p}_i}$ is now the vector formed by column stacking that matrix. 
Similarly, $\mathbf{Y}_{\mathbf{c}_i}$ and $\mathbf{y}_{\mathbf{c}_i}$ are, respectively, the matrix and column-stacked vector containing measurements at spatial channel $i$, $\mathbf{c}_i$, at all parameters and all times. 
Then $\mathbf{Y}_{t_i}$ and $\mathbf{y}_{t_i}$ are the matrix and column-stacked vector containing all spatial channel measurements for all parameters at time $t_i$.
When constructing the questionnaire metric for the parameters, for example, the main change is in the definition of the basis vectors, which now must be defined over the spatial channel indices and the time channel indices. To do this, we cluster the spatial channels \{$\mathbf{c}_i$\} into a set of clusters  \{$\mathbf{I}_l$\} and the temporal channels \{$t_j$\} into another set of clusters \{$\mathbf{J}_{l^\prime}$\}. Then for all $l$,$l^\prime$, we define

\begin{equation}
\mathbf{g}_{l,l^\prime}(\mathbf{c}_i,t_j)=
\left\{
    \begin{array}{ll}
        1, & \mathbf{c}_i \in \mathbf{I}_l, t_j \in \mathbf{J}_{l^\prime} \\
        0, & otherwise
    \end{array}
\right.
\end{equation}

\noindent Instead of summing over a cluster, our basis vectors now sum over the intersection between a spatial channel cluster and a temporal channel cluster. Our iterative algorithm becomes:
\begin{enumerate}[nosep]
	\item Use an uninformed metric to cluster the parameter channels and the spatial channels
	\item Use the uninformed parameter and channel clusters to define an informed metric and clustering for the time channels.
	\item Repeat until converged:
	\begin{enumerate}[nosep]
	    \item Use the spatial channel and time channel clusters to update the informed parameter channel metric and clusters
	    \item Use the parameter channel and time channel clusters to update the spatial channel metric and clusters
	    \item Use the parameter channel and spatial channel clusters to update the time channel metric and clusters
    \end{enumerate}
\end{enumerate}
\noindent  Once again, this is a brief operational summary (adapted to our example) of the general questionnaire approach introduced in Ref. \citenum{Ankenman2014} and used in, e.g., Ref. \citenum{Yair2017,Sroczynski2018,Mishne2016}.

\section{\textit{DROSOPHILA} EMBRYO MODEL}
\label{app:drosophila}

\begin{figure}[h!]
	\begin{center}
		\includegraphics[width=8.7cm]{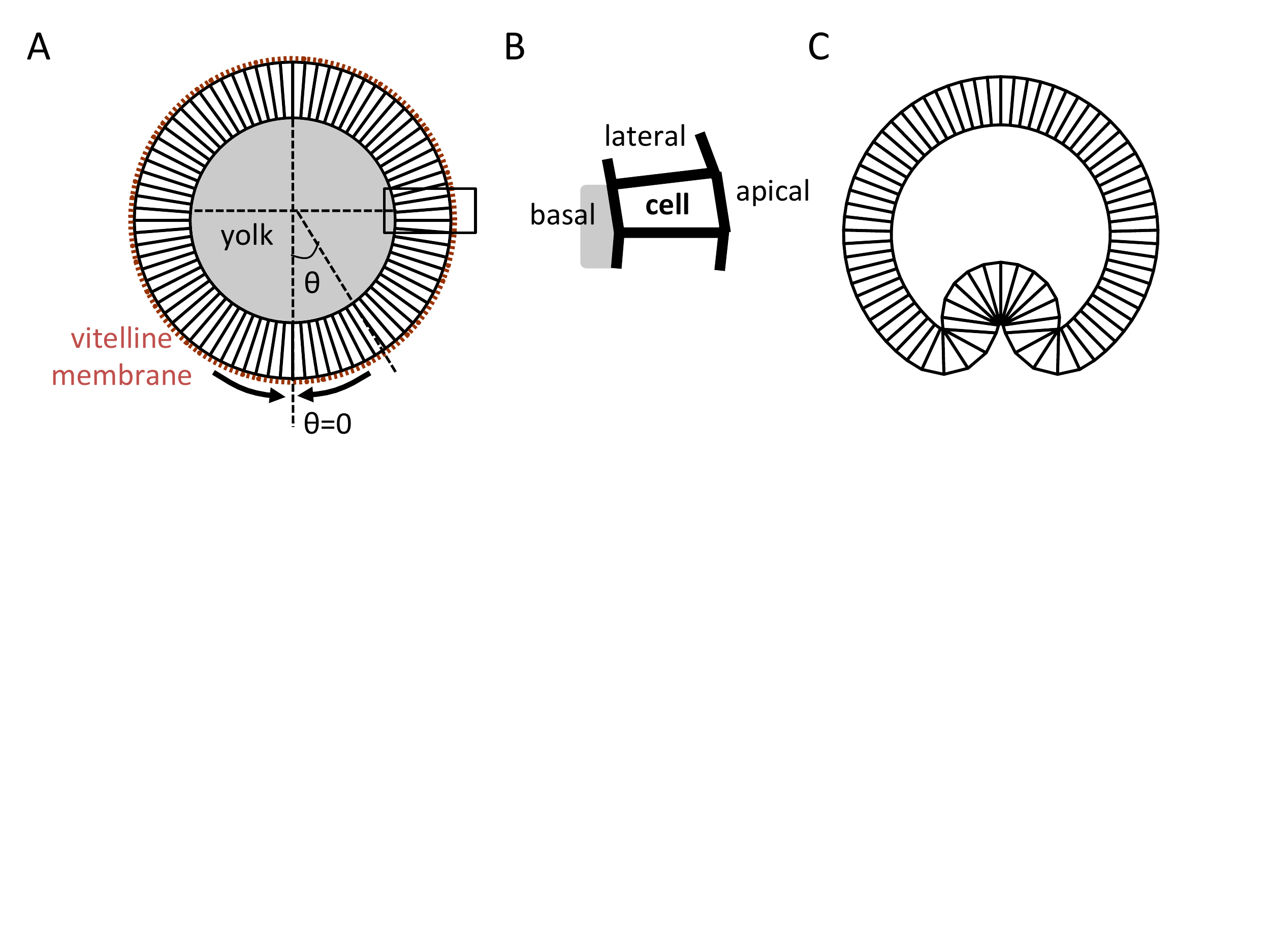}
		\caption{Schematic representation of the 2D {\em Drosophila} egg cross-section vertex model. (A) Initial configuration consisting of a ring of $N_c$ identical trapezoidal cells of area $A_c^0$, enclosing a yolk of area $A_Y^0$, and surrounded by a membrane. (B) A schematic highlighting a single cell, with distinct apical, basal and lateral edges. (C) A representative configuration after introduction of additional patterning of model parameters, which resembles an embryo cross-section after ventral furrow formation.}
		\label{f1}
	\end{center}
\end{figure}

\begin{figure*}
    \includegraphics[width=11.4cm]{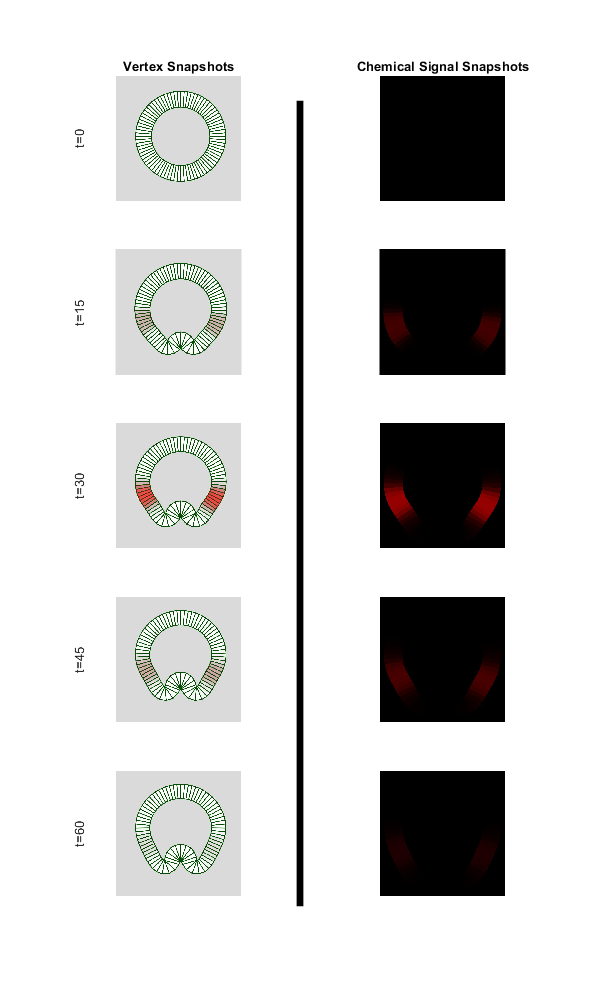}
    \caption{\label{fig:example_movie} \textit{(Left)} Representative temporal snapshots showing the evolution of the location of cell vertices as well as of the chemical signal intensity in the cells. \textit{(Right)} Corresponding observed snapshots (the observer can only see the chemical signal from chemical staining).}
\end{figure*}

\subsection{Energy formulation of the computational model}
In Ref. \citenum{Misra2016}, the authors modeled the evolution of the $Drosophila$ embryo cross-section near the midpoint of its anterior-posterior axis, at the onset of ventral furrow formation, as a two-dimensional cross-section of a cylindrical monolayered epithelium (Fig. \ref{f1} A). The cells in the embryo are modeled as quadrilaterals with distinct apical and basal edges; each lateral edge is shared between two cells. The ring of cells encloses the yolk and is surrounded by a stiff membrane.

Guided by previous studies \cite{HocevarBrezavscek2012,Polyakov2014,Misra2016,Misra2016a,Rauzi2013, Odell1981}, the following energy formulation was used to model the epithelium:

\begin{equation}
E_{2D}=\sum_a\sigma_al_a+\sum_b\sigma_bl_b+\sum_l\sigma_ll_l+B\sum_c(A_c-A_c^0)^2
\end{equation}

In this expression, the first three terms sum over all apical edges $a$, basal edges $b$, and lateral edges $l$, respectively. These terms capture intercellular interactions in the form of tensile forces; $l_a,l_b,$ and $l_l$ are the edge lengths of the apical, basal, and lateral edges, respectively, and $\sigma_a,\sigma_b,$ and $\sigma_l$ are the corresponding line tension coefficients. The last term sums over all cells $c$ and penalizes the deviation of the cell area $A_c$ from its target value $A_c^0$. $B$ is the compression modulus approximating cytoplasmic incompressibility.

The effects of the outer stiff membrane and the enclosed yolk are captured by introducing additional terms. The total energy of the system is given by the following expression:

\begin{equation}
E = B_Y(A_Y-A_Y^0)^2 + \epsilon\sum_k\frac{1}{(R_c-R_k)^n} + E_{2D}
\end{equation}

The first term penalizes deviations from the initial area of the inner cavity; $A_Y$ is the area enclosed by the epithelium, $A_Y^0$ is the initial area, and $B_Y$ is the compression modulus. The second term sums over all apical vertices $k$ and models the outer stiff membrane by restricting the radial motion of the vertices within a circle of radius $R_c$ that is concentric with the initial homogeneous configuration. $R_k$ denotes the radial distance of apical vertex $k$, and hence, $R_c-R_k$ represents the membrane thickness near that vertex. The value $\epsilon$ is the membrane stiffness parameter, and $n$ is the exponent of the repulsive potential.

The epithelial monolayer is assumed to evolve in an over-damped setting. The dynamics of the vertex $j$ with position vector ${\bf{x}}_j$ are governed by the following equation:

\begin{equation}
\eta \frac{d{\bf{x}}_j}{dt}= {\bf{F}}_j= -\nabla_jE
\end{equation}                                                

$\eta$ is the mobility coefficient and $F_j$ is the force acting on vertex $j$. The equations are propagated in time using direct forward-Euler integration implemented in C++.

\subsection{Geometric parameters of the initial homogeneous configuration}

All simulations start from an initial homogeneous configuration consisting of $N_c=80$ trapezoidal cells with identical geometry and model parameters. A spatial pattern of cell properties is later imposed to drive shape transformation. To ensure that this configuration is mechanically stable,  the energy of a single cell is minimized in an idealized setting.  $h$ is defined as the height of the trapezoid, which is assumed to be approximately equal to $l_l$. By taking the cell area compression modulus to be large ($B\rightarrow\infty$), one can assume that $A_c$ remains similar to $A_c^0$, providing a relationship between geometric variables $l_a,l_b,$ and $h$ ($h=\frac{2A_c^0}{l_a+l_b}$). The equation for a single cell reduces to 

\begin{equation}
E_c=\sigma_l\frac{2A_C^0}{l_a+l_b}+\sigma_al_a+\sigma_bl_b
\end{equation}

The minima of this function ($\frac{\partial E_c}{\partial l_a} = 0$, $\frac{\partial E_c}{\partial l_b}=0$, and positive definite Hessian matrix) give a relationship between the material properties and the geometric parameters for the stable homogeneous configuration.

\subsection{Shape transformation driven by modifications to apical and basal line tension}

A position dependent value of the apical line tension coefficient $\sigma_a$ is imposed for the apical edges located on the ventral side. The value of the coefficient $\sigma_a$ for each apical edge depends on the angle the mid-point of the edge makes from the reference line (i.e., $\theta=0$, Fig. \ref{f1} A). Hence, the value of $\sigma_a$ depends on the position of the edges and can change with time as the configuration evolves.

The spatial pattern of $\sigma_a$ is given by the following function:
\begin{gather}
\sigma_a =
\left\{
	\begin{array}{ll}
		\sigma_{a,0}(1+P\cdot e^{-\theta^2/G^2})  & \mbox{if } |\theta| < \pi/4\\
		\sigma_{a,0} & \mbox{if } |\theta| \geq \pi/4
	\end{array}
\right.
\end{gather}

$\sigma_{a,0}$ corresponds to the apical line tension coefficient for the stable homogeneous configuration. The value of $\sigma_a$ is peaked at the ventral-most point and the spatial pattern is present within an angular spread of ventral cells.

A reduction of the basal line tension coefficient for all basal edges ($\sigma_b=f\cdot\sigma_{b,0}$, where $\sigma_{b,0}$ correspond to the parameter value for a stable homogeneous configuration and $f<1$) was found necessary to simulate a more closed form of the ventral furrow. 

\subsection{Chemical signal model}

The concentration $C$ of the chemical signal is governed by the following partial differential equation

\begin{equation}
\frac{\partial C}{\partial t} = D\frac{\partial^2 C}{\partial x^2} + r(t)G(x) -dC
\end{equation}

Here, $r(t)$ represent the rate of production, $G(x)$ is a binary-valued function (1 if the region in space is a source of the signal, 0 otherwise), $D$ is the diffusion coefficient, and $d$ is the degradation coefficient.

The following form was chosen for the rate of production:
\begin{equation}
r =
\left\{
	\begin{array}{ll}
		kt^2  & \mbox{if } 0 \leq t < t_s \\
		kt_s^2e^{-\alpha(t-t_s)} & \mbox{if } t \geq t_s
	\end{array}
\right.
\end{equation}

The rate of production $r$ increases polynomially till $t=t_s$ and then decreases exponentially with coefficient $\alpha$.

Space was discretized into 80 bins corresponding to the 80 cells and approximate the diffusion operator using a second order finite difference scheme with periodic boundary conditions: 

\begin{equation}
\frac{\partial C_i}{\partial t} = D_e(C_{i+1}-2C_i+C_{i-1})+r(t)G(i)-dC_i,
\end{equation}

where $i=1,2,\ldots 80$ is the cell index and $D_e$ is an effective diffusion coefficient. Indexing of the cells is done such that the ventral most cells are indexed as 1 and 80 (lying on opposide sides of the reference line $\theta=0$). The above system of ordinary differential equations was solved using the forward-Euler method.

\subsection{Plotting morphology and signal}

\begin{figure}[h!]
	\begin{center}
		\includegraphics[scale=0.5]{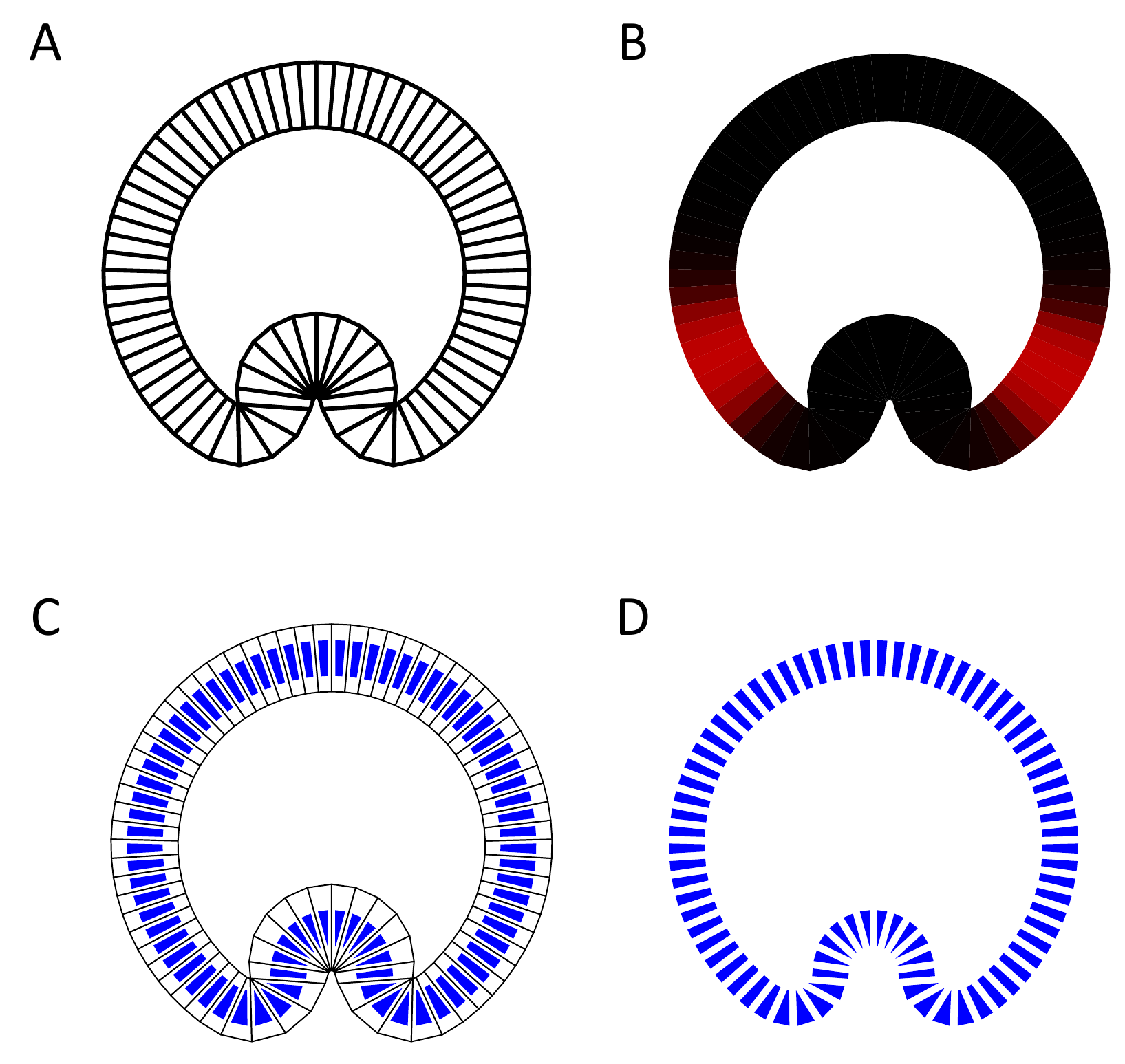}
		\caption{ Morphology and the chemical signal. (A) A representative configuration during the simulation. (B) A representative plot of the chemical signal intensity. (C) An overlay of the exact model configuration and the cell {\it{nuclei}} (shrunk versions of actual cells). (D) A representative plot of the morphology displaying only the nuclei. }
		\label{f2}
	\end{center}
\end{figure}

The model state at a given time is defined by the coordinates of each vertex (${\textbf{x}}_j$) and the concentration of signal in each cell ($C_i$). To create visualizations that approximate experimental images,  the concentration was scaled based on the maximum concentration that the signal reaches across a set of simulations and fill in each cell red with intensity proportional to the scaled concentration (Fig. \ref{f2} B). Since experimental images contain signals corresponding to the nucleus of each cell, one plots the nucleus by shrinking each cell to a fraction of its actual area and filling in the resulting quadrilateral blue (Fig. \ref{f2} C,D). 

\subsection{Parameter values}

Geometric parameters of intial configuration:

Number of cells $N_c$ = 80

Lateral edge length $l_l$ = 3.0

Radius of apical circle $R_a$ is a normally distributed variable with mean = 8.5 and standard deviation = 0.25 (with a constraint that it lies within 8.25 and 8.75)

The above geometric variable fixes other geometric parameters.

Radius of the basal circle $R_b=R_a-l_l$, apical edge length $l_a = (2\pi R_a)/N_c$, and basal edge length $l_b = (2\pi R_b)/N_c$.

Model parameters for initial homogeneous configuration:

$\sigma_{a,0}=2.6418,\sigma_{b,0}=2.6418,\sigma_l=1.0,B=20.0,B_Y=0.01,\epsilon=10^{-10},n=4,R_c=1.05*R_a$

Model parameters for additional patterning:

$P=0.2, G=4.0, f=0.7$

Parameters used for simulations:

$D_e=0.2, k=5*10^{-5}, \alpha =0.03, t_s=40,$

$d$ is a normally distributed variable with mean = 0.08 and standard deviation = 0.008 (with a constraint that it lies within 0.056 and 0.104).

$G(i)=1$ for $i=13-19,61-67$, otherwise $G(i)=0$.

\section{Learning a non-autonomous PDE}\label{app:non-autonomous}

\begin{figure*}
    \includegraphics[width=\textwidth,trim={5cm 9cm 5cm 0},clip]{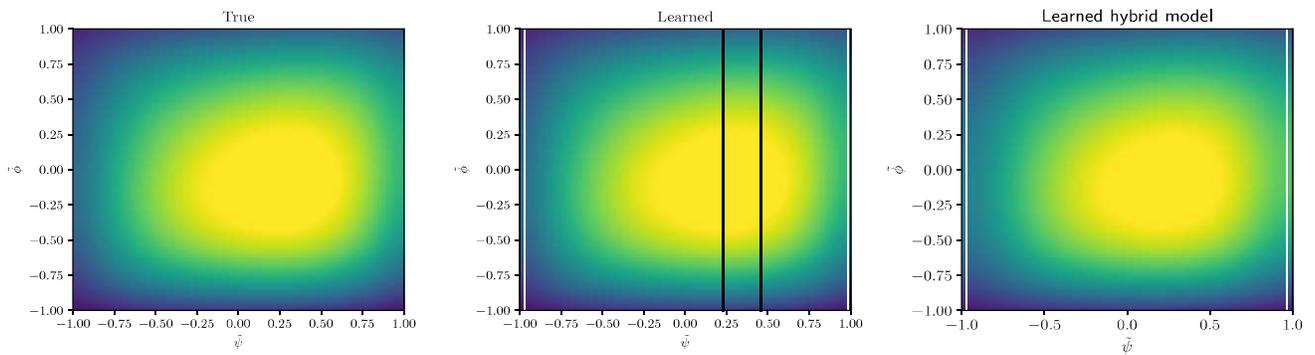}
    \caption{Left: True data as shown in Fig.~\ref{fig:learn_time}(right). Center: Data obtained by integrating the same initial snapshot of the data as shown on the left, but using the learned model $f$. Boundary ``corridor'' conditions are provided, as indicated by the white vertical lines. In addition, a corridor around the source term is provided, as indicated by the black vertical lines. Right: Integration results using the hybrid (partly autonomous, partly temporally forced) model, ~\eqref{eq:f_w_source}, using the same initial snapshot. Note the visual agreement with the true data,}
    \label{fig:learn_results_emergent_space_and_time_plus_hybrid}
\end{figure*}

In the main body of the paper, we learned the autonomous part of the dynamical system {\em outside of the source corridor}, and provided the true source corridor data to act as boundary conditions for our integrations. Alternatively, one can extend this approach {\em to learn the source term {\underline in the corridor} $I_s$} as well
\begin{equation}
    source~term = g\left(\left.\tilde{\psi}\right|_{\tilde{\psi} \in I_s}, \tilde{\phi}\right),
\end{equation}
such that the predictions of the learned autonomous PDE model~\eqref{eq:tau_model}, are augmented to be accurate in the source corridor $I_s$. That is, keeping the learned (outside the source corridor) model $f$ fixed, we can learn the source term of a non-autonomous dynamical system
\begin{equation}
\frac{\partial c}{\partial \tilde{\phi}} = f\left(c, \frac{\partial c}{\partial \tilde{\psi}}, \dots, \frac{\partial^n c}{\partial \tilde{\psi}^n}\right) + g\left(\left.\tilde{\psi}\right|_{\tilde{\psi} \in I_s}, \tilde{\phi}\right)
\label{eq:f_w_source}
\end{equation}
with $c=c\left(\tilde{\psi}, \tilde{\phi}\right)$ now using training data inside the source corridor.
In particular, we represent $g$ using a fully connected neural network with three hidden layers of 64 neurons each, each followed by a Swish activation function and using the hyper parameters described above.
The weights are optimized using the mean squared error between the true $\partial c/\partial \tilde{\phi}$ and the predictions of the fixed autonomous PDE $f$ in the source interval $I_s$.
This means that after training, $g$ maps from a point in emergent space in the source interval, $I_s$, and time $\tau$ to the difference $\partial c/\partial \tilde{\phi} - f$ at that point in emergent space and emergent time.
Note that $g\left(\left.\tilde{\psi}\right|_{\tilde{\psi} \in I_s}, \tilde{\phi}\right)=0$ outside the source interval as marked by the black vertical lines in Fig.~\ref{fig:learn_results_emergent_space_and_time}.
Integrating this time dependent model, Eq.~\ref{eq:f_w_source}, the predictions are indeed accurate over the
entire interval considered (see integration results in Fig.~\ref{fig:learn_results_emergent_space_and_time_plus_hybrid}(right)). This means that we learned a source term that is consistent with the data and the learned autonomous spatial operator over the entire domain. 
For the machine-learned PDE to work, some important qualitative features (lack of translational invariance, non-autonomousness) {\em must be qualitatively known to be present in advance}, so that that the corresponding functional forms can be correctly identified. 

\end{document}